\documentclass[12pt]{article}

\usepackage{amsfonts,latexsym,amssymb,amsmath}

\def \al{\alpha}
\def \be{\beta}

\def \ep{\varepsilon}

\def \ka{\varkappa}

\def \ph{\varphi}

\def \Om{\Omega}

\def \operatorname#1{\mathop{\rm #1}}

\def\div{\operatorname{div}}
\def\osc{\operatorname{osc}}
\def\osc2{\operatorname{osc^2}}

\def\cd{\partial}

\def\Q0{Q(x_0,t_0,R)}
\def\0{{x_0,t_0,R}}
\def\build#1_#2{\mathrel{\mathop{\kern 0pt#1}\limits_{#2}}}

\newtheorem{theorem}{Theorem}[section]

\newtheorem{proposition}{Proposition}[section]
\newtheorem{lemma}{Lemma}[section]

\title{On the Stokes Problem  with  \\ Non-Zero Divergence
\thanks{This work is supported by
RFBR grant 08-01-00372-a.}}
\author{N.~Filonov, T.~Shilkin}
\date{}
\begin{document}

\maketitle

\centerline{\it Dedicated to Nina Nikolaevna Uraltseva}

\abstract{We study the strong solvability of the nonstationary Stokes problem with non-zero divergence in a bounded domain. }

\section{Introduction and Main Results}
\setcounter{equation}{0}

Let $\Om$ be a domain in $\mathbb R^n $, $n\ge 2$,
with sufficiently smooth boundary $\cd\Om$,
and assume that $\Om$ is homeomorphic to a ball.
We study the solvability of the linear initial boundary-value problem
\begin{equation}
\left. \begin{array}c
\cd_t v -\Delta v +\nabla p = f \\ \div v =g
\end{array}
\right\}\qquad \mbox{in } \ Q_T:=\Om\times(0,T)\\
\label{Problem_1}
\end{equation}
\begin{equation}
v|_{t=0}=0,\qquad v|_{\cd\Om\times(0,T)}=0.
\label{IC_BC}
\end{equation}
We assume there are $s$, $l\in (1,+\infty)$ such that the following conditions hold:
\begin{equation}
f\in L_{s, l}(Q_T), \label{External force}
\end{equation}
\begin{equation}
g \in W^{1,0}_{s,l}(Q_T),  \label{Assumptions_1}
\end{equation}
\begin{equation}
\cd_t g \in L_{s,l}(Q_T), \label{Assumptions_2}
\end{equation}
\begin{equation}
  \int\limits_{\Om} g(x,t)~dx =0 ,\qquad \mbox{a.e. } t\in(0,T), \qquad g(\cdot,0)=0.
  \label{Assumptions_3}
\end{equation}
Here  $L_{s,l}(Q_T)$ is the anisotropic Lebesgue space equipped with the norm
$$
\|f\|_{L_{s,l}(Q_T)}:=
\Big(\int_0^T\Big(\int_\Om |f(x,t)|^s~dx\Big)^{l/s}dt\Big)^{1/l} ,
$$
and we use the following notation for the functional spaces:
$$
\gathered
W^{1,0}_{s,l}(Q_T)\equiv L_l(0,T; W^1_s(\Om))= \{ \ u\in L_{s,l}(Q_T): ~\nabla u \in L_{s,l}(Q_T) \ \},
\\
W^{2,1}_{s,l}(Q_T) = \{ \ u\in W^{1,0}_{s,l}(Q_T): ~\nabla^2 u, \ \cd_t u \in L_{s,l}(Q_T) \ \},
\\
\overset{\circ}{W}{^1_s}(\Om)=\{ \ u\in W^1_s(\Om):~ u|_{\cd\Om}=0 \ \},
\\
W^{-1}_s(\Om)=(\overset{\circ}{W}{^1_{s'}}(\Om))^*= \ \mbox{ dual space to } \overset{\circ}{W}{^1_{s'}}(\Om),
\endgathered
$$
and the following notation for the norms:
$$
\gathered
\| u \|_{W^{1,0}_{s,l}(Q_T)}= \| u \|_{L_{s,l}(Q_T)}+ \|\nabla u\|_{L_{s,l}(Q_T)},  \\
\| u \|_{W^{2,1}_{s,l}(Q_T)}= \| u \|_{W^{1,0}_{s,l}(Q_T)}+ \| \nabla^2 u \|_{L_{s,l}(Q_T)}+\|\cd_t u\|_{L_{s,l}(Q_T)},  \\
\| u\|_{W^{-1}_s(\Om)} =  \sup\limits_{\overset{ w\in \overset{\circ}{W}{^1_{s'}}(\Om)  }{ \|\nabla w\|_{L_{s'}(\Om)}\le 1 } }\Big|\int_{\Om} u\cdot w~dx\Big|,  \\
\| u \|_{L_l(0,T; W^{-1}_s(\Om))} = \Big(\int_0^T \| u(\cdot,t) \|^l_{W^{-1}_s(\Om)}~dt\Big)^{1/l}.
\endgathered
$$
Our main result is the following
\begin{theorem} \label{Nonzero divergence}
Assume $s$, $l\in (1,\infty)$ and let  $f$,  $g$ satisfy conditions  (\ref{External force}) --- (\ref{Assumptions_3}).
 Then there exists the unique pair of functions $(v,\nabla p)$ such that
 \begin{equation*}
 v\in W^{2,1}_{s,l}(Q_T),\qquad \nabla p\in L_{s,l}(Q_T),
 \end{equation*}
and $(v,\nabla p)$  satisfy the equations (\ref{Problem_1}) a.e. in  $Q_T$ and  (\ref{IC_BC}) in the sense of traces. Moreover, the following estimate holds:
\begin{equation}
\gathered
\| v \|_{W^{2,1}_{s,l} (Q_T)} + \| \nabla
p\|_{L_{s,l}(Q_T)} \le \\ \le C_*\left(\|f\|_{L_{s,l}(Q_T)}+ \|  g  \|_{W^{1,0}_{s,l}(Q_T)} +
\| \cd_t g\|_{L_{s,l}(Q_T)}^{1/s}  \|
\cd_t g \|_{L_{l}(0, T; W^{-1}_s(\Om))}^{1/s'}\right).
\endgathered
\label{Estimate_1}
\end{equation}
Here $C_*$ is a constant depending only on   $n$, $T$, and  \  $\Om$.

\end{theorem}

The following theorem shows that the assumption \eqref{Assumptions_2}
in Theorem \ref{Nonzero divergence}
can not be omitted or replaced  by a weaker assumption
 \begin{equation}
 \cd_t g\in L_{l}(0,T; W^{-1}_s(\Om)). \label{Assumptions_4}
  \end{equation}

\begin{theorem}\label{Negative result}
Assume $n=2$ and $\Om$ is a unit disc in $\mathbb R^2$.
There exist  functions $f$, $g$ satisfying conditions \eqref{External force}, \eqref{Assumptions_1}, \eqref{Assumptions_3}, \eqref{Assumptions_4} with $s=l=2$
and $g|_{\cd\Om \times (-1, 0)} = 0$,
and there exists a weak solution  $(v, p)$ of the problem \eqref{Problem_1}
in $Q=\Om\times (-1,0)$ satisfying  the initial data $v|_{t=-1}=0$
and the boundary data $v|_{\cd\Om\times (-1,0)}=0$ in the sense of traces,
and possessing the properties
\begin{equation}
v\in C([0,T]; L_2(\Om))\cap  W^{1,0}_2(Q_T),
\label{Class_Weak}
\end{equation}
\begin{equation}
p\in L_2(Q_T),
\label{Pressure}
\end{equation}
\begin{equation}
\cd_t v \in L_2(0,T; W^{-1}_2(\Om)),
\label{Time Derivative}
\end{equation}
\begin{equation}
(v,p,f,g) \quad \mbox{satisfy  (\ref{Problem_1}) in the sense of distributions,}
\label{Traces}
\end{equation}
but
$$
v\not \in W^{2,1}_2(Q), \qquad\nabla p \not\in L_2(Q) ,
$$
so the weak solution $(v,p)$ fails to be a strong solution.
\end{theorem}

Theorem \ref{Negative result} exhibits  nonexistence of a strong solution
to the problem (\ref{Problem_1}), (\ref{IC_BC}) under the assumptions
(\ref{External force}), (\ref{Assumptions_1}), (\ref{Assumptions_3}) (\ref{Assumptions_4}) only, as the following uniqueness theorem shows:

\begin{theorem} \label{Uniqueness of weak solution}
Assume $n\ge 2$ and $f$, $g$ satisfy conditions (\ref{External force}), (\ref{Assumptions_1}), (\ref{Assumptions_3}), (\ref{Assumptions_4}) with $s=l=2$. Then the weak solution of the problem (\ref{Problem_1}), (\ref{IC_BC}) possessing the properties (\ref{Class_Weak})---(\ref{Traces}) (if exists) is unique.
\end{theorem}

The counterexample provided by Theorem \ref{Negative result}
looks surprising as if we take an arbitrary divergent-free
function $v$ such that
$$
v\in W^{2,1}_2(Q_T), \quad v|_{\cd\Om}=0, \quad  v|_{t=0}=0,
$$
then we have
$$
\cd_t \div v \ \in \ L_2(0,T; W^{-1}_2(\Om)) ,
$$
and one could conjecture that condition \eqref{Assumptions_4}
with $l=s=2$ is the natural one for the solvability of the problem
\eqref{Problem_1}, \eqref{IC_BC} in the class
$(v,p)\in W^{2,1}_2(Q_T)\times W^{1,0}_2(Q_T)$.
Theorem \ref{Negative result} demonstrates that this is not the case.

\medskip
Estimates of  Sobolev norms of a solution $v$ to the  problem (\ref{Problem_1}) by  Lebesgue  norms of the functions $f$, $\nabla g$ and $\cd_t g$  are well-known, see,
for example, \cite{Farwig_Sohr}.
The specific feature of our estimate (\ref{Estimate_1}) is its multiplicative form,
i.e.   right-hand side of (\ref{Estimate_1}) includes a product of a stronger norm
$\| \cd_t g\|_{L_{s,l}(Q_T)}$ by a weaker norm
$\|\cd_t g \|_{L_{l}(0, T; W^{-1}_s(\Om))}$.
Such form is convenient for a simple proof of the local estimates
of solutions of the Stokes problem near the boundary:

\begin{proposition}\label{Stokes_Local}
Denote $Q^+:=\{ x\in \mathbb R^n :  |x|<1 ,~ x_n>0 \} \times (-1,0)$
and
$$
Q^+_{1/2} :=\{ x\in \mathbb R^n : |x|<1/2 ,~ x_n>0 \} \times (-1/4,0).
$$
Assume $u\in W^{2,1}_{s,l}(Q^+)$, $q\in W^{1,0}_{s,l}(Q^+)$,
$\tilde f\in L_{s,l}(Q^+)$ satisfy the following Stokes system:
\begin{equation}
\gathered
\left.
\begin{array}c
\cd_t u -\Delta u +\nabla q = \tilde f  \\
\div u = 0
\end{array}\quad
\right\} \qquad\mbox{ in } \ Q^+,  \\
u|_{x_n=0}=0.
\endgathered
\label{Stokes_near_boundary}
\end{equation}
Then there is an absolute constant $C$ (depending only on $n$) such that
\begin{equation}
\gathered
\| u \|_{W^{2,1}_{s,l}(Q^+_{1/2})}+\| \nabla q\|_{L_{s,l}(Q^+_{1/2})}\le
\\ \le C \Big(\|\tilde f\|_{L_{s,l}(Q^+)}+
\|u \|_{W^{1,0}_{s,l}(Q^+)}+ \inf_{b\in L_l(-1,0)} \| q-b\|_{L_{s,l}(Q^+)}\Big) .
\label{Local_Stokes}
\endgathered
\end{equation}
\end{proposition}

We remark that estimate \eqref{Local_Stokes} plays an important role in
the study of the boundary regularity of suitable weak solutions to
the Navier-Stokes system, see \cite{Seregin_JMFM}, \cite{Seregin_Handbook}
and reference there.
The estimate (\ref{Local_Stokes}) was proved  in \cite{Seregin_ZNS_271}. In \cite{Solonnikov_ZNS_288} the same result was established for the generalized Stokes system. The local Stokes problem (\ref{Stokes_near_boundary}) can be transferred to the initial boundary-value problem of type (\ref{Problem_1}) by multiplication of $u$ by appropriate cut-off function $\zeta$, where  $v=\zeta u$, $p=\zeta q$.
Then the estimate (\ref{Local_Stokes}) follows easily from (\ref{Estimate_1}) by
iterations.
We reproduce the derivation of (\ref{Local_Stokes}) from  (\ref{Estimate_1}) in the Appendix of the present paper.

Theorem \ref{Nonzero divergence} gives only sufficient conditions
for the solvability of the problem \eqref{Problem_1}
in the class $W^{2,1}_{s,l}(Q_T)$.
The  conditions on $g$ which are both necessary and sufficient for
the strong solvability of the problem \eqref{Problem_1}
seems to be  unknown even in the case of $s=l=2$.

In \cite{Solonnikov_Uspekhi} the following estimate was proved for solution $(v,p)$ of the problem (\ref{Problem_1}), (\ref{IC_BC}):
\begin{equation}
\gathered
\| v \|_{W^{2,1}_{s,l} (Q_T)} + \| \nabla
p\|_{L_{s,l}(Q_T)} \le \\ \le  C_*\left(\| f\|_{L_{s,l}(Q_T)}+\| \nabla g  \|_{L_{s,l}(Q_T)} +
 \| \cd_t g \|_{L_{l}(0, T; \hat W^{-1}_s(\Om))}\right),
\endgathered
\label{Estimate_2}
\end{equation}
where $\| \cdot\|_{\hat W^{-1}_s(\Om)}$ stands for the dual norm to the space $W^1_{s'}(\Om)$ (with non-zero traces on the boundary):
\begin{equation*}
 \| v\|_{\hat W^{-1}_s(\Om)} =  \sup\limits_{\overset{ w\in {W}{^1_{s'}}(\Om)  }{ \| w\|_{W^1_{s'}(\Om)}\le 1 } }\Big|\int_{\Om} v\cdot w~dx\Big| .
\end{equation*}
We remark that the estimate (\ref{Estimate_2}) is not so convenient
for applications as  a weak solution $u\in W^{1,0}_{s,l}(Q^+)$,
$q\in L_{s,l}(Q^+)$ of the local Stokes problem (\ref{Stokes_near_boundary})
satisfies the estimate
\begin{equation}
\| \cd_t u \|_{ L_l(-1,0; W^{-1}_s(B^+))}\le
C (\|\tilde f\|_{ L_l(-1,0; W^{-1}_s(B^+))}
+\|u \|_{W^{1,0}_{s,l}(Q^+)}+\|q \|_{L_{s,l}(Q^+)} )
\label{119}
\end{equation}
but, generally speaking, the similar estimate with $\|\cd_t u \|_{ L_l(0,T;  W^{-1}_s(B^+))}$ replaced by  $\|\cd_t u \|_{ L_l(0,T; \hat W^{-1}_s(B^+))}$  is not true.

Our paper is organized as follows:
in Section \ref{Axillary Results} we present several auxiliary theorems
concerning extensions of functions from the boundary onto a whole domain;
in Section \ref{On the problem div u = g} we prove a theorem on solutions
to the problem $\div u =g$, $u|_{\cd\Om}=0$;
the proof of Theorem \ref{Nonzero divergence} is presented in
the Section \ref{Section 2};
a counterexample of Theorem \ref{Negative result} is constructed
in Section \ref{Section 3};
in the Appendix the derivation of the estimate \eqref{Local_Stokes}
from (\ref{Estimate_1}) is given.

\section{Auxiliary Results}\label{Axillary Results}
\setcounter{equation}{0}

\bigskip
In this section we formulate several results concerning extension theorems
from the boundary of a domain.
We denote by $\mathbb R^n_+$ the half-space
$\mathbb R^n_+ = \{ x = (x', x_n) \in \mathbb R^n : x_n>0\}$,
and by $\nabla '$ the gradient with respect to $x'$.
Let us start with the following

\begin{proposition} \label{Boundary_multiplicative}
For any $\ph\in W^1_s(\Om)$ the following estimate holds:
\begin{equation}
\| \ph\|_{L_s(\cd\Om)}\le C \| \ph\|_{L_s(\Om)}^{1/s'}\| \ph\|_{W^1_s(\Om)}^{1/s}.
\label{BM_estimate}
\end{equation}
\end{proposition}

\medskip\noindent{\bf Proof.}
For a function $\ph:\mathbb R^n_+\to \mathbb R$
the estimate \eqref{BM_estimate}
follows from the integral representation
$$
|\ph(x',0)|^s=-\int\limits_0^{+\infty} \frac{\cd}{\cd x_n}|\ph(x',x_n)|^s~dx_n
$$
with the help  of the H\" older inequality.
For a bounded smooth domain $\Om\subset \mathbb R^n$ the estimate \eqref{BM_estimate} can be justified by a standard techniques
of the local maps and partition of unity.
$\blacksquare$

By $W^r_s (\cd\Om)$ with non-integer $r>0$ we denote
the Sobolev-Slo\-bo\-det\-skii
space of functions defined on $\cd\Om$.
The next proposition is essentially  proved in \cite{Sol_Ur}.
We just need to verify that the extension operator $T_1$
can be constructed in such a way that both estimates
\eqref{Extension_1} and \eqref{Extension_2} hold simultaneously.

\begin{proposition} \label{Extension proposition 2}
Let $\Om\subset \mathbb R^n$ be a bounded domain, $\cd\Om \in C^3$.
There exists a linear operator $T_1$
 $$
 T_1: W^{2-\frac 1s}_s(\cd\Om)\times W^{1-\frac 1s}_s(\cd\Om) \to W^2_s(\Om)
 $$
 such that for any $ b\in W^{2-\frac 1s}_s(\cd\Om)$,
$ a\in W^{1-\frac 1s}_s(\cd\Om)$ the function $f:=T_1(b,a)$ possesses the following properties:
\begin{equation*}
\gathered
f|_{\cd\Om} = b, \qquad
\frac{\cd f}{\cd \nu}\Big|_{\cd\Om}=a,
\endgathered
\end{equation*}
\begin{equation}
\| f \|_{W^1_s(\Om)}\le C_1 \Big(\| b \|_{W^1_s(\cd\Om)} + \| a \|_{L_s(\cd\Om)}\Big).
\label{Extension_1}
\end{equation}
Moreover, if additionally  $b\in W^{3-\frac 1s}_s(\cd\Om)$,
$a\in W^{2-\frac 1s}_s(\cd\Om)$ then $f \in W^3_s(\Om)$ and
\begin{equation}
\| f \|_{W^3_s(\Om)}\le C_2 \Big(\| b \|_{W^{3-\frac 1s}_s(\cd\Om)} + \| a \|_{W^{2-\frac 1s}_s(\cd\Om)}\Big).
\label{Extension_2}
\end{equation}
The constants $C_1$ and $C_2$ depend only on $n$ and $\Om$.
\end{proposition}

\medskip\noindent{\bf Proof.}
First, we consider the case of a half-space, $\Omega = \mathbb R^n_+$.
Assume $a\in W^{1-\frac 1s}_s(\mathbb R^{n-1})$ and
$b\in W^{2-\frac 1s}_s(\mathbb R^{n-1})$.
Let us consider a kernel  $K\in C_0^\infty(\mathbb R^{n-1})$ with the following properties:
$$
\int\limits_{\mathbb R^{n-1}} K(y')~dy'=1,
\qquad \int\limits_{\mathbb R^{n-1}} y_\al K(y')~dy'=0,
\quad\al = 1, \ldots, n-1,
$$
and a smooth cut-off function $\zeta: [0,+\infty)\to \mathbb R$ such that
$$
\zeta(y_n)\equiv 1 \quad\mbox{on}\quad [0,1/2], \quad 0\le \zeta\le 1, \qquad \zeta(y_n)\equiv 0 \quad\mbox{on}\quad [1,+\infty).
$$
Define the function $f$ as follows:
\begin{gather*}
f (y) = \zeta(y_n) (g (y)+ h (y)),\\
g (y) = \int\limits_{\mathbb R^{n-1}} K(z') b (y'+y_nz')~dz',\\
h (y) = y_n\int\limits_{\mathbb R^{n-1}} K(z') a (y'+y_n z')~dz'.
\end{gather*}
Then obviously $f|_{y_n=0} = b$,
$\frac{\cd f}{\cd y_n} |_{y_n=0} = a$.
It is well known that for
$a\in W^{2-\frac 1s}_s(\cd\Om)$,
$b\in W^{3-\frac 1s}_s(\cd\Om)$,
the inequality
\begin{equation*}
\| f \|_{W^3_s(\Om)}\le C_2 \Big(\| b \|_{W^{3-\frac 1s}_s(\cd\Om)} + \| a \|_{W^{2-\frac 1s}_s(\cd\Om)}\Big)
\end{equation*}
holds (see \cite{Sol_Ur}).
So, we need to verify the estimate
\begin{equation}
\| f \|_{W^1_s(\mathbb R^n_+)}
\le C\Big( \| b \|_{W^1_s(\mathbb R^{n-1})}
+ \| a \|_{L_s(\mathbb R^{n-1})}\Big).
\label{Weaker estimate}
\end{equation}
Consider, for example, the function $h$.
We have
\begin{gather*}
h(y) = y_n^{2-n} \int\limits_{\mathbb R^{n-1}}
K \left(\frac{z'-y'}{y_n}\right) a (z')~dz', \\
\frac{\cd h (y)}{\cd y_\alpha} = y_n^{1-n} \int\limits_{\mathbb R^{n-1}}
\frac{\cd K}{\cd y_\alpha} \left(\frac{z'-y'}{y_n}\right) a (z')~dz', \\
\frac{\cd h (y)}{\cd y_n} = y_n^{1-n} \int\limits_{\mathbb R^{n-1}}
\left( (2-n) K \left(\frac{z'-y'}{y_n}\right) -
\langle \nabla ' K \left(\frac{z'-y'}{y_n}\right) , \frac{z'-y'}{y_n} \rangle \right)
a (z')~dz' .
\end{gather*}
Integral convolution operators in $L_s$-spaces are bounded
by $L_1$-norm of the kernel.
Therefore,
\begin{gather*}
\| \zeta h \|_{L_s(\mathbb R^n_+)}
\le \|K\|_{L_1(\mathbb R^{n-1})} \| a \|_{L_s(\mathbb R^{n-1})} , \\
\left\| \frac{\cd (\zeta h)}{\cd y_\alpha} \right\|_{L_s(\mathbb R^n_+)}
\le \left\|\frac{\cd K} {\cd y_\alpha} \right\|_{L_1(\mathbb R^{n-1})}
\| a \|_{L_s(\mathbb R^{n-1})} , \quad \al = 1, \ldots, n-1, \\
\left\| \frac{\cd (\zeta h)}{\cd y_n} \right\|_{L_s(\mathbb R^n_+)}
\le C \| a \|_{L_s(\mathbb R^{n-1})} \\
\end{gather*}
and $\| \zeta h \|_{W^1_s(\mathbb R^n_+)}
\le C \| a \|_{L_s(\mathbb R^{n-1})} $,
where the constant $C$ can be explicitly expressed in terms of functions
$K$ and $\zeta$.
The inequality $\| \zeta g \|_{W^1_s(\mathbb R^n_+)}
\le C \| b \|_{W^1_s(\mathbb R^{n-1})}$ follows by the similar argument.
Thus, we justified \eqref{Weaker estimate}.

Again, the case of a bounded smooth domain reduces to
the case of a half-space by the standard techniques of localisation.
$\blacksquare$

\bigskip
Now we formulate one result from \cite{Solonnikov_PMA_21}.
This result is an analog of Bogovskii's result \cite{Bogovskii} in the case of smooth compact manifold $\cd\Om$.
Assume $\Om\subset \mathbb R^n$ is a domain which is homeomorphic to a ball
and denote by $\nu(x)$ the unit outer normal to $\cd\Om$ at the point $x\in \cd\Om$.
Let $b:\cd\Om\to \mathbb R^n$ be a vector field such that $b\cdot \nu=0$.
Below the symbol ${\rm div}_{S}~ b $ stands for the differential operator
which is defined in a local coordinate system $\{y_\al\}_{\al=1}^{n-1}$  by
$$
{\rm div}_{S}~ b = \frac 1{\sqrt{g}}\frac{\cd}{\cd y_\al} (\sqrt{g}~ \hat b_\al(y)),
$$
where $g=\det (g_{\al\be})$, \
$g_{\al\be}=\frac{\cd x(y)}{\cd y_\al}\cdot \frac{\cd x(y)}{\cd y_\be}$,
and $\hat b_\al(y)$ are the components of a vector field $b$
in local coordinates $\{ y_\al\}$,
i.e. $b(x(y)) = \hat b_\al(y)\frac {\cd x(y)}{\cd y_\al}$.

\begin{proposition} \label{Surface div}
Assume $\Om\subset \mathbb R^n$ is a smooth domain which is homeomorphic to a ball.  There exists a linear operator $T_2$
 $$
 T_2:\{ \ \ka\in W^{1-\frac 1s}_s(\cd\Om): \ \int_{\cd\Om} \ka ~ ds =0   \ \}\to W^{2-\frac 1s}_s(\cd\Om;\mathbb R^n),
 $$
 such that
the function \ $b = T_2 \ka $ \ possesses the following properties:
$$
\langle b, \nu\rangle = 0 , \qquad
{\rm div}_{S}~ b= \ka\quad \mbox{on}\quad \cd\Om,
$$
and
\begin{equation}
\| b \|_{W^{1}_s(\cd\Om)}\le C \| \ka \|_{L_s(\cd\Om)}.
\label{b weaker}
\end{equation}
Moreover, if additionally  $\ka \in W^{2-\frac 1s}_s(\cd\Om)$ then
\begin{equation}
\| b \|_{W^{3-\frac 1s}_s(\cd\Om)}\le C \| \ka \|_{W^{2-\frac 1s}_s(\cd\Om)}.
\label{b stronger}
\end{equation}
\end{proposition}

\noindent
Proposition \ref{Surface div} is proved in \cite{Solonnikov_PMA_21},
see Propositions 2.1, 2.2, 2.3 there.
We just emphasize that as the construction of the operator $T_2$
in a local coordinates $\{ y_\al\}$ uses nothing but the Bogovskii operator
(see \cite{Bogovskii}), the both estimates (\ref{b weaker}) and (\ref{b stronger})
are satisfied simultaneously.

Combining Propositions \ref{Extension proposition 2}
and \ref{Surface div}, we finally obtain

\begin{proposition}\label{rot proposition}
Let $\Om\subset \mathbb R^n$ be a bounded domain
which is homeomorphic to a ball, $\cd\Omega \in C^4$.
Then there exists a linear operator
$$
T_3: \{ \ \ka\in W^{1-\frac 1s}_s(\cd\Om):
\ \int_{\cd\Om}\ka ~ds=0 \ \} \to W^2_s(\Om;\mathbb R^n),
$$
such that the function $w =T_3\ka$ possesses the properties
\begin{equation*}
\div w = 0, \quad w|_{\cd\Om} = -\ka \, \nu , \quad
\| w \|_{L_s(\Om)} \le C \| \ka \|_{L_s(\cd\Om)}.
\end{equation*}
Moreover, if additionally  $\ka\in W^{2-\frac 1s}_s(\cd\Om)$ then
$ w \in W^2_s(\Om;\mathbb R^n)$ and
$\| w \|_{W^2_s(\Om)} \le C \| \ka \|_{W^{2-\frac 1s}_s(\cd\Om)}$.
\end{proposition}

\medskip\noindent{\bf Proof.}
Denote by $\tilde \nu$ a smooth extension of the field $\nu$
into the whole domain $\Omega$, $\tilde \nu : \Omega \to \mathbb R^n$,
$\tilde \nu |_{\cd\Omega} = \nu$.
Let
$$
b=-T_2\ka \in W^{2-\frac 1s}_s(\cd\Om; \mathbb R^n),
\quad \langle b, \nu\rangle =0 .
$$
Define the vector-field
$$
a = \langle b, \nabla \rangle \tilde \nu
- b \div \tilde \nu \in W^{2-\frac 1s}_s(\cd\Om),
$$
and let $f = T_1 (b, a)$, where $T_1$ is the operator constructed in Proposition
\ref{Extension proposition 2}.
We have
\begin{equation}
f|_{\cd\Om} = b, \qquad
\frac{\cd f}{\cd \nu}\Big|_{\cd\Om}=a,
\label{p1}
\end{equation}
\begin{equation}
\| f \|_{W^1_s(\Om)} \le C \| b \|_{W^1_s(\cd\Om)}
\le \tilde C \| \ka \|_{L_s(\cd\Om)},
\label{p2}
\end{equation}
and
\begin{equation}
\| f \|_{W^3_s(\Om)} \le C \| b \|_{W^{3-\frac 1s}_s(\cd\Om)}
\le \tilde C \| \ka \|_{W^{2-\frac 1s}_s(\cd\Om)}
\label{p3}
\end{equation}
in the case $\ka \in W^{2-\frac 1s}_s(\cd\Om)$.
Note that $|\nu(x)|^2 = 1$ on the boundary,
so $\langle b, \nabla \rangle \nu \perp \nu$ and
$\frac{\cd f}{\cd\nu} = a \perp \nu$ on $\cd\Omega$.
Therefore,
\begin{equation}
(\div f)|_{\cd\Om} = {\rm div}_{S}~ b .
\label{p4}
\end{equation}
Now we introduce the vector-function
$w\in W_s^1 (\Omega, \mathbb R^n)$ defined as
$$
w_j (x)= \sum_{i=1}^n \frac{\cd}{\cd x_i}
\left( f_i(x) \tilde\nu_j(x) - f_j(x) \tilde\nu_i(x) \right) .
$$
Clearly, $\div w = 0$.
We have also
$\| w \|_{L_s(\Om)} \le C \| \ka \|_{L_s(\cd\Om)}$ and
$$
\| w \|_{W^2_s(\Om)} \le C \| \ka \|_{W^{2-\frac 1s}_s(\cd\Om)},
\quad \ka\in W^{2-\frac 1s}_s(\cd\Om),
$$
due to \eqref{p2} and \eqref{p3}.
Finally, by virtue of \eqref{p1} and \eqref{p4} we get
\begin{gather*}
w|_{\cd\Om} = \left. \left( \tilde \nu \div f + \langle f, \nabla \rangle \tilde \nu
- \langle \tilde \nu, \nabla \rangle f - f \div \tilde \nu \right) \right|_{\cd\Om} \\
= \nu \, {\rm div}_{S}~ b + a - \frac{\cd f}{\cd \nu} = - \nu \ka .
\qquad\blacksquare
\end{gather*}

\section{On the  problem $\div u=g$} \label{On the problem div u = g}
\setcounter{equation}{0}

\begin{theorem} \label{div u=g}
There exists a linear operator
$$
T: \{\ g \in L_s(\Om) : \ \int_\Om g\, dx=0 \ \} \quad \to \quad \overset{\circ}{W}{^1_s}(\Om; \mathbb R^n)
$$
such that the function $u=Tg$ is a solution of the equations
$$
\left\{ \begin{array}c
\div u = g \quad\mbox{a.e. in}\quad \Om \\
u|_{\cd\Om}=0
\end{array}
\right.
$$
which satisfies the estimate
\begin{equation*}
\| u \|_{L_s(\Om)}\le C_1 \| g\|_{L_s(\Om)}^{1/s}\| g \|_{W^{-1}_s(\Om)}^{1/s'}.
\end{equation*}
Moreover, if $g\in W^1_s(\Om)$ then $u\in W^2_s(\Om)$ and
$\| u \|_{W^2_s(\Om)}\le C_2 \| g \|_{W^1_s(\Om)}$.
Here $C_1$ and $C_2$ depend only on $n$, $s$, and $\Om$.
\end{theorem}

\noindent{\bf Proof.}
Let $\ph\in \overset{\circ}{W}{^1_s}(\Om)\cap W^2_s(\Om)$
be a solution to the Dirichlet problem
\begin{equation*}
\Delta \ph =g \quad\mbox{in}\quad \Om, \qquad \ph|_{\cd\Om}=0,
\end{equation*}
and define the function $\ka:\cd\Om\to\mathbb R$ by the formula
$\ka = \frac{\cd\ph}{\cd\nu}$.
We have
$$
\| \ph\|_{W^1_s(\Om)}\le C \| g\|_{W^{-1}_s(\Om)}, \quad
\| \ph\|_{W^2_s(\Om)}\le C \| g\|_{L_s(\Om)}
$$
and by Proposition \ref{Boundary_multiplicative}
$\| \ka \|_{L_s(\cd\Om)}\le C
\| g\|_{L_s(\Om)}^{1/s}\| g \|_{W^{-1}_s(\Om)}^{1/s'}$.
If $g\in W^1_s(\Om)$ then
$$
\| \ka \|_{W^{2-\frac 1s}_s(\cd\Om)}
\le C \| \ph\|_{W^3_s(\Om)} \le \tilde C \| g\|_{W^1_s(\Om)}.
$$
Note that
$\int_{\cd\Om} \ka ds = \int_\Om g dx = 0$,
so we can apply Proposition \ref{rot proposition} to the function $\ka$.
Let $w = T_3 \ka$ and $u = \nabla \ph + w$.
Then
$$
\|u\|_{L_s(\Om)} \le \|\ph\|_{W_s^1(\Om)} + C \|\ka\|_{L_s(\cd\Om)}
\le C_1 \| g\|_{L_s(\Om)}^{1/s}\| g \|_{W^{-1}_s(\Om)}^{1/s'}
$$
and
$$
\|u\|_{W^2_s(\Om)} \le \|\ph\|_{W_s^3(\Om)} + C \|\ka\|_{W_s^{2-\frac 1s}(\cd\Om)}
\le C_2 \| g \|_{W^1_s(\Om)}.
$$
Finally,
$u|_{\cd\Om} = \frac{\cd\ph}{\cd\nu} \nu - \ka \nu = 0$.
$\blacksquare$

\section{Proof of Theorem \ref{Nonzero divergence}}\label{Section 2}
\setcounter{equation}{0}

Assume $g$ satisfies conditions \eqref{Assumptions_1} -- \eqref{Assumptions_3}
and consider the function $w=Tg$,
where the operator $T$ is defined in Theorem \ref{div u=g}.
Then
\begin{gather*}
\div w =g  \quad \mbox{a.e. in}\quad Q_T, \qquad w|_{\cd\Om \times (0,T)}=0, \\
w(\cdot, 0) =0,\qquad \cd_t w= T(\cd_t g) \quad \mbox{a.e. in}\quad Q_T, \\
\| w (\cdot, t)\|_{W^2_s(\Om)}\le
C \| g(\cdot, t)\|_{W^1_s(\Om)}\qquad \mbox{for a.e. } t\in (0,T), \\
\| \cd_t w(\cdot, t)\|_{L_s(\Om)}\le
C \| \cd_t g(\cdot, t)\|_{L_s(\Om)}^{1/s}
\| \cd_t g(\cdot, t)\|_{W^{-1}_s(\Om)}^{1/s'}\quad
\mbox{for a.e. } t\in (0,T).
\end{gather*}
Taking the power $l$, integrating these inequalities with respect to $t$
and applying the H\" older inequality, we obtain
\begin{equation}
\| w\|_{W^{2,1}_{s,l}(Q_T)}\le C \Big(\|g\|_{W^{1,0}_{s,l}(Q_T)}+ \| \cd_t g\|_{L_{s,l}(Q_T)}^{1/s}\| \cd_t g\|_{L_l(0,T; W^{-1}_s(\Om))}^{1/s'}\Big) .
\label{function w}
\end{equation}
Let  $(u,\nabla p)$ be the solution of the  Stokes problem
\begin{equation*}
\gathered
\left.\begin{array}c
\cd_t u - \Delta u +\nabla p = f - (\cd_t w -\Delta w) \\
\div u = 0
\end{array}\right\}\qquad\mbox{in}\quad Q_T
\\
u|_{\cd \Om\times(0,T) } = 0, \qquad  u|_{t=0}=0.
\endgathered
\end{equation*}
It is well-known (see, for example, \cite{Solonnikov_Uspekhi} and references there)
that $(u,\nabla p) $ satisfy the estimate
\begin{equation}
\| u \|_{W^{2,1}_{s,l}(Q_T)}+\| \nabla p\|_{L_{s,l}(Q_T)}\le C\Big(\| f\|_{L_{s,l}(Q_T)} + \|w\|_{W^{2,1}_{s,l}(Q_T)} \Big) .
\label{Stokes_homogeneous}
\end{equation}
Put $v=u+w$.
Then $(v,\nabla p)$ is a solution to the problem \eqref{Problem_1},
\eqref{IC_BC}.
Combining estimates \eqref{function w} and \eqref{Stokes_homogeneous}
we obtain \eqref{Estimate_1}.
$\blacksquare$

\section{Proofs of Theorems \ref{Negative result} and
\ref{Uniqueness of weak solution}}\label{Section 3}
\setcounter{equation}{0}

For the  presentation convenience in this section we denote by $\Om$
the unit disc in $\mathbb R^2$ and by $Q\subset \mathbb R^2\times \mathbb R$
we denote the following space-time cylinder
$$
Q:=\Om\times (-1,0).
$$
Moreover, we assume the Stokes system (\ref{Problem_1}) is considered in $Q$
and the initial value $v|_{t=-1}=0$ is prescribed at $t=-1$.

\bigskip\noindent
{\bf Proof of Theorem \ref{Negative result}.}

\noindent
{\bf 1.}
For $t<0$ we introduce the scalar function $\psi:Q\to \mathbb R$
given by serie
$$
\psi(r,\theta, t) := \sum\limits_{n=1}^\infty
\frac{r^n \sin n\theta}{n^4 (1-n^7t)}
$$
in the polar coordinate system
$x_1 = r \cos\theta$, $x_2 = r\sin\theta$.
Then
$$
\cd_r \psi(r,\theta)= \sum\limits_{n=1}^\infty \frac{r^{n-1} \sin n\theta}{n^3 (1-n^7t)}, \qquad
\frac 1r \cd_\theta \psi =\sum\limits_{n=1}^\infty \frac{r^{n-1} \cos n\theta}{n^3 (1-n^7t)}
$$
and $\Delta \psi = 0$ in $Q$.
Introduce the vector-function $w:Q\to \mathbb R^2$ which is given
by formulas \ $\vec w=w_r \vec e_r + w_\theta \vec e_\theta$,
$$
w_r(r,\theta,t):=
\sum\limits_{n=1}^\infty \frac{\al_n(r) \sin n\theta}{n^3 (1-n^7t)}, \qquad
w_\theta(r,\theta,t):=
\sum\limits_{n=1}^\infty \frac{\al_n(r) \cos n\theta}{n^3 (1-n^7t)}.
$$
Here $\al_n\in W^2_\infty(0,1)$ are any functions satisfying the following  conditions:
\begin{equation}
\al_n(r) =\left\{ \begin{array}{cl} 0, & r\in [0,1-\frac 1{n^3}],  \\
0< \al_n(r)< 1, & r\in (1-\frac 1{n^3}, 1), \\
\al_n(r) = 1, & r=1.
\end{array}\right.
\label{Conditions_for_alpha_n_1}
\end{equation}
\begin{equation}
\al'_n(1) = n-1,\\\
\label{Conditions_for_alpha_n_2}
\end{equation}
\begin{equation}
 |\al'_n(r)|\le Cn^3, \quad |\al''_n(r)|\le Cn^6\qquad\forall~r\in [0,1].
\label{Conditions_for_alpha_n_3}
\end{equation}
For example, the following functions $\al_n$ satisfy all conditions (\ref{Conditions_for_alpha_n_1}) --- (\ref{Conditions_for_alpha_n_3}):
$$
\alpha_n (r) = (3n^6 - n^4 + n^3) (r - 1 + n^{-3})^2 - (2n^9 - n^7 + n^6) (r - 1 + n^{-3})^3
$$
for \ $r\in (1-\frac 1{n^3}, 1]$ \ and \ $\al_n(r)=0$ \ for \ $r\in  [0,1-\frac 1{n^3}]$.

Take a smooth cut-off function in $t$-variable $\chi\in C^1([-1,0])$  such that
$$
0\le \chi(t)\le 1,
\quad
\chi(t)=0 \quad \forall ~t\in [-1,- 2/3 ],
\quad
\chi(t) =1 \quad \forall~ t\in [- 1/3, 0],
$$
and denote by $v$, $p$, $f$,  $g$ the following functions:
\begin{equation}
\gathered
v:= \chi (w-\nabla \psi), \qquad p:=\chi \cd_t \psi, \\
f: = \chi (\cd_t w-\Delta w) +\chi' (w-\nabla\psi),
\qquad g:=\chi \div w.
\endgathered
\label{formulas for u,p,f,g}
\end{equation}
Then $(v,p,f,g)$ satisfy pointwise the following system of equations:
\begin{equation}
\gathered
\left. \begin{array}c
\cd_t v -\Delta v +\nabla p = f \\ \div v =g
\end{array}
\right\}\qquad \mbox{in } \ Q=\Om\times(-1,0)\\
\label{Problem_2}
v|_{t=-1}=0,\qquad v|_{\cd\Om}=0.
\endgathered
\end{equation}
Moreover, for any $t\in (-1,0)$ we have
$$
\int_\Om g(x,t)~dx=\chi(t) \int_{\cd\Om}w(s,t)\cdot \nu(s)~ds= \chi(t)\int_0^{2\pi}w_r(1,\theta, t)~d\theta=0.
$$
>From (\ref{Conditions_for_alpha_n_2}) we obtain
$$
\gathered
\div w~\big|_{\cd\Om} = \Big(\cd_r w_r +\frac 1r w_r +\frac 1r\cd_\theta w_\theta\Big) ~\Big|_{r=1}=  \\ =
\sum\limits_{n=1}^\infty \left.\frac{(\al_n'+\frac {\al_n}r-n\al_n) \sin n\theta}{n^3 (1-n^7t)}~\right|_{r=1} =0
\endgathered
$$
So, $g|_{\cd\Om \times (-1, 0)} = 0$.

\bigskip
\noindent
{\bf 2.} Below we will  show that the following relations hold:
\begin{equation}
\chi w\in W^{2,1}_2(Q),
\label{to be verified 1}
\end{equation}
\begin{equation}
\chi\psi \in W^{2,1}_2(Q),
\label{to be verified 4}
\end{equation}
\begin{equation}
\cd_t\nabla (\chi \psi) \not\in L_2(Q).
\label{to be verified 6}
\end{equation}
These relations imply that the data $(f,g)$ of the problem \eqref{Problem_2}
given by formulas \eqref{formulas for u,p,f,g} possess all the properties \eqref{Class_Weak} -- \eqref{Traces}.
But this weak solution is not a strong one as
$\cd_t v \not\in L_2(Q)$ and $\nabla p\not\in L_2(Q)$.

We start from the verification of \eqref{to be verified 1}.
We have
$$
\cd_t w_r = \sum\limits_{n=1}^\infty \frac{n^4 \al_n(r) \sin n\theta}{ (1-n^7t)^2}
$$
and hence
$$
\gathered
\| \cd_t w_r\|_{L_2(Q)}^2=  \int\limits_{-1}^0 dt \int\limits_{0}^{2\pi} d\theta \int\limits_{0}^{1} |\cd_t w_r(r,\theta)|^2~rdr =
\pi \sum\limits_{n=1}^\infty \int\limits_{-1}^0
\int\limits_0^1 \frac{n^8 |\al_n(r)|^2 rdrdt}{ (1-n^7t)^4}
\endgathered
$$
As $\int\limits_0^1 |\al_n(r)|^2~ rdr\le n^{-3}$
we obtain
 $$
\gathered
\| \cd_t w_r\|_{L_2(Q)}^2 \le  C \sum\limits_{n=1}^\infty \ \int\limits_{-1}^0  \ \frac{n^5 dt}{ (1-n^7t)^4}\le C \sum\limits_{n=1}^\infty \frac 1{n^2}<+\infty.
\endgathered
$$
A similar estimate holds for $\| \cd_t w_\theta\|_{L_2(Q)}$.
Hence we conclude $\cd_t w \in L_2(Q)$.
Now we turn to the estimate of $\| \nabla^2 w\|_{L_2(Q)}$:
$$
\gathered
\| \nabla^2 w\|_{L_2(Q)}^2 \le
C \sum\limits_{n=1}^\infty \ \int\limits_{-1}^0 \int\limits_0^1 \
\frac{ \left( |\al_n''|^2 r + n^2 |\al_n'|^2 r^{-1} + n^4 |\al_n|^2 r^{-3} \right)
~ drdt}{ n^6 (1-n^7 t)^2} .
\endgathered
$$
The conditions \eqref{Conditions_for_alpha_n_1} and \eqref{Conditions_for_alpha_n_3}
imply
$$
\int_0^1
\left( |\al_n''|^2 r + n^2 |\al_n'|^2 r^{-1} + n^4 |\al_n|^2 r^{-3} \right) ~ dr
\le C n^9 ,
$$
so
$$
\gathered
\| \nabla^2 w\|_{L_2(Q)}^2\le C \sum\limits_{n=1}^\infty \
\int\limits_{-1}^0  \ \frac{ \ n^3~dt}{  (1-n^7t)^2}\le  C \sum\limits_{n=1}^\infty \frac {1}{n^4}<+\infty.
\endgathered
$$
The weaker norms $\|w\|_{L_2(Q)}$ and $\|\nabla w\|_{L_2(Q)}$
can be estimated in the similar way.
So, \eqref{to be verified 1} is proved.
The proof of \eqref{to be verified 4} is analogous.

We are left to prove (\ref{to be verified 6}).
>From  (\ref{to be verified 4}) we see that $\chi'\nabla \psi\in L_2(Q)$
and hence we need to show that $\chi \cd_t\nabla\psi\not\in L_2(Q)$.
As $\chi \equiv 1$ on $[-\frac 13,0]$
and the functions $\{ \sin n\theta \}_{n=1}^\infty$ are orthogonal in $L_2(0,2\pi)$
it is sufficient to show that
\begin{equation}
\sum\limits_{n=1}^\infty \ \int\limits_{-1/3}^0~dt \ \int\limits_{0}^1 \
\left( \frac{n^4 r^{n-1}}{(1 - n^7 t)^2} \right)^2 ~ rdr = +\infty .
\label{Nonconvergence}
\end{equation}
Indeed,
$$
\int\limits_{-1/3}^0~dt \ \int\limits_{0}^1 \
\frac{n^8 r^{2n-1} dr}{(1 - n^7 t)^4}
= \frac{1}{6} + O (n^{-21}), \quad n \to \infty,
$$
thus we arrive at \eqref{Nonconvergence}.
$\blacksquare$

\bigskip
\noindent
{\bf Proof of Theorem \ref{Uniqueness of weak solution}.}
Assume there are two week solutions $(v_1, p_1)$ and $(v_2,p_2)$
satisfying the system \eqref{Problem_1}, \eqref{IC_BC}
with the same functions $(f,g)$.
Consider the differences $w=v_1-v_2$, $q=p_1-p_2$.
Then $(w,q)$ is a weak solution to the homogeneous Stokes problem with zero data.
This solution satisfies all conditions \eqref{Class_Weak}---\eqref{Traces}.
Multiplying the equation by $w$ we obtain
$$
\frac 12 \, \cd_t \|w\|_{L_2}^2 = - \|\nabla w \|_{L_2}^2 \le 0,
$$
and therefore $w\equiv 0$.
$\blacksquare$

\section{Appendix}
\setcounter{equation}{0}

In this section we present the derivation of the estimate \eqref{Local_Stokes}
from the estimate \eqref{Estimate_1}.
We remind that $Q^+:=B^+\times (-1,0)$,
$B^+:=\{ \ x\in\mathbb R^n:~|x|<1, \ x_n>0 \ \}$
and take arbitrary $\rho$, $r$ such that
$$
\begin{array}c
\frac12\le\rho<r\le\frac{9}{10}.
\end{array}
$$
Consider a cut-off function $\zeta\in C_0^\infty(Q)$ such that
$$
\gathered
0\le \zeta\le 1 \quad\mbox{in} \quad Q^+,
\qquad\zeta\equiv 1 \quad\mbox{in} \quad Q^+_\rho,
\qquad \zeta\equiv 0 \quad\mbox{in} \quad Q^+\setminus Q^+_r, \\
 \| \nabla^k\zeta \|_{L_\infty(Q^+)}\le \frac {C}{(r-\rho)^{k}},\quad k=1,2,
\quad \| \cd_t\zeta \|_{L_\infty(Q^+)}\le \frac {C}{r-\rho},
\endgathered
$$
where
$$
Q^+_R:= B^+_R\times (-R^2,0), \quad
B^+_R:= \{ x\in \mathbb R^n: |x|<R, x_n > 0 \} .
$$
Let $(u,q)$ be a solution to the system (\ref{Stokes_near_boundary}) and consider functions
$v:=\zeta u$, $p:=\zeta q$.
Then $(v,p)$ is a solution to the problem (\ref{Problem_1})
with $\Om $ being a smooth domain such that \
$B^+_{{9}/{10}}\subset \Om\subset B^+_1$ \  and
$$
f = \zeta \tilde f+ u (\cd_t \zeta- \Delta \zeta) - 2(\nabla u)\nabla \zeta+ q\nabla \zeta , \quad g = u \cdot \nabla \zeta .
$$
Applying the estimate (\ref{Estimate_1}) and taking into account that
$\frac 1{r-\rho}\ge 1$ we obtain
$$
\gathered
\| u \|_{W_{s,l}^{2,1} (Q_\rho^+)}^s \le C\| \tilde f\|_{L_{s,l}(Q^+)}^s +
\frac C{(r-\rho)^{2s}} \Big(\| u \|_{W_{s,l}^{1,0}(Q^+)}^s +
\|q\|_{L_{s,l}(Q^+)}^s\Big) + \\
+C\Big( \| \nabla (u\cdot\nabla \zeta) \|_{L_{s,l}(Q^+)}^s + \| \cd_t (u\cdot\nabla \zeta) \|_{L_{s,l}(Q^+)}\| \cd_t (u\cdot\nabla \zeta) \|_{L_l(-1,0; W^{-1}_s(B^+))}^{s-1}\Big) .
\endgathered
$$
Taking into account estimates
$$
\gathered
\| \nabla (u\cdot\nabla \zeta) \|_{L_{s,l}(Q^+)}^s\le \frac C{(r-\rho)^{2s}}\|u \|_{W^{1,0}_{s,l}(Q^+)}^s, \\
\| \cd_t (u\cdot\nabla \zeta) \|_{L_{s,l}(Q^+)}\le
\frac C{(r-\rho)^2} \Big(\| \cd_t u \|_{L_{s,l}(Q^+_r)}
+ \| u \|_{L_{s,l}(Q^+)}\Big), \\
\| \cd_t (u\cdot\nabla \zeta) \|_{L_l(-1,0; W^{-1}_s(B^+))}^{s-1}\le \frac C{(r-\rho)^{2s-2}} \Big(\| \cd_t u \|_{L_l(-1,0; W^{-1}_s(B^+))}^{s-1} + \|  u \|_{L_{s,l}(Q^+)}^{s-1}\Big),
\endgathered
$$
we get
\begin{equation}
\gathered
\| u \|_{W_{s,l}^{2,1} (Q_\rho^+)}^s
\le C\| \tilde f\|_{L_{s,l}(Q^+)}^s  \\
+ \frac C{(r-\rho)^{2s}} \Big(\| u \|_{W^{1,0}_{s,l}(Q^+)}^s
+ \|q\|_{L_{s,l}(Q^+)}^s+ \| \cd_t u \|_{L_l(-1,0; W^{-1}_s(B^+))}^{s}\Big) \\
+ \frac C{(r-\rho)^{2s}} \| \cd_t u \|_{L_{s,l}(Q^+_r)}\Big(
\| \cd_t u \|_{L_l(-1,0; W^{-1}_s(B^+))}^{s-1}+ \|  u \|_{L_{s,l}(Q^+)}^{s-1}\Big) .
\endgathered
\label{To_chto_nado}
\end{equation}
Estimating the last term in the right-hand side of \eqref{To_chto_nado}
via the Young inequality $ab\le \ep a^s+C_\ep b^{s'}$ we obtain the estimate
$$
\gathered
\frac C{(r-\rho)^{2s}} \| \cd_t u \|_{L_{s,l}(Q^+_r)}\Big(
\| \cd_t u \|_{L_l(-1,0; W^{-1}_s(B^+))}^{s-1}+ \|  u \|_{L_{s,l}(Q^+)}^{s-1}\Big) \le \\ \le \ep \| \cd_t u \|_{L_{s,l}(Q^+_r)}^s + \frac{C_\ep}{(r-\rho)^{2ss'}}\Big(
\| \cd_t u \|_{L_l(-1,0; W^{-1}_s(B^+))}^{s}+ \|  u \|_{L_{s,l}(Q^+)}^{s}\Big) ,
\endgathered
$$
where the constant $\ep>0$ can be chosen arbitrary small.
Therefore,
\begin{gather*}
\| u \|_{W_{s,l}^{2,1} (Q_\rho^+)}^s
\le  C\| \tilde f\|^s_{L_{s,l}(Q^+)}+ \ep \| \cd_t u \|_{L_{s,l}(Q^+_r)}^s + \\
 +
\frac {C_\ep}{(r-\rho)^{2ss'}} \Big(\| u \|_{W^{1,0}_{s,l}(Q^+)}^s  + \|q\|_{L_{s,l}(Q^+)}^s+ \| \cd_t u \|_{L_l(-1,0; W^{-1}_s(B^+))}^{s}\Big) ,
\end{gather*}
and by virtue of \eqref{119}
\begin{equation}
\gathered
\| u \|_{W_{s,l}^{2,1} (Q_\rho^+)}^s
\le  \ep \| \cd_t u \|_{L_{s,l}(Q^+_r)}^s  \\
 +
\frac {C_\ep}{(r-\rho)^{2ss'}}
\Big(\| \tilde f\|_{L_{s,l}(Q^+)}^s + \| u \|_{W^{1,0}_{s,l}(Q^+)}^s + \|q\|_{L_{s,l}(Q^+)}^s \Big) .
\endgathered
\label{inequality this implies}
\end{equation}
Now let us introduce the monotone function
$\Psi(\rho) := \| u \|_{W_{s,l}^{2,1} (Q_\rho^+)}^s$,
and the constant
$$
A:=C_\ep\left(\| \tilde f\|_{L_{s,l}(Q^+)}^s + \| u \|_{W^{1,0}_{s,l}(Q^+)}^s  + \|q\|_{L_{s,l}(Q^+)}^s\right).
$$
The inequality \eqref{inequality this implies} implies that
\begin{equation}
\begin{array}c
\Psi (\rho) \le \ep\Psi(r)+\frac {A} {(r-\rho)^\al}, \qquad \forall ~\rho, \ r:
\quad R_1 \le \rho<r\le R_0,
\end{array}
\label{Giaquinta's lemma}
\end{equation}
for some  $\al>0 $ depending only on $s$, and for $R_1=\frac 12$,  $R_0=\frac 9{10}$.
Now we shall take an advantage of  the following lemma (which can be easily  proved by iterations if one take $r_k:=R_0-2^{-k}(R_0-R_1)$):

\begin{lemma}
Assume $\Psi$ is a nondecreasing bounded function which satisfies the inequality (\ref{Giaquinta's lemma}) for some $\al>0$, $A>0$, and $\ep\in (0,2^{-\al})$.
Then there exists a constant $B$ depending only on  $\ep$ and $\al$ such that
$$
\Psi(R_1)\le \frac {B\, A}{(R_0-R_1)^\al} .
$$
\label{Giaquinta_lemma}
\end{lemma}
Fixing $\ep = 2^{-3 ss'}$ in \eqref{inequality this implies} and applying
Lemma \ref{Giaquinta_lemma} to our function $\Psi$,
we obtain the estimate
$$
\| u \|_{W_{s,l}^{2,1} (Q_{1/2}^+)} \le
C_* \Big(\| \tilde f\|_{L_{s,l}(Q^+)}^s+\| u \|_{W^{1,0}_{s,l}(Q^+)}^s  + \|q\|_{L_{s,l}(Q^+)}^s \Big)
$$
which completes the proof.
$\blacksquare$


\begin{thebibliography}{99}

\bibitem{Besov}
{\sc O.V.~Besov, V.P.~Il'in, S.M.~Nikolskii},  Integral
representations of functions and imbedding theorems. Nauka,
Moscow, 1975. Translation: Wiley\&Sons, 1978.

\bibitem{Bogovskii} {\sc M.E.~Bogovskii}, {\it On solution of some
problems of vectoral analysis related to  $\div$ and $\operatorname{grad}$
operators}, Proc. of S.L.~Sobolev Seminar {\bf 1} (1980), 5-40.

\bibitem{CKN}
{\sc L.~Caffarelli, R.V.~Kohn, L.~Nirenberg}, {\it Partial
regularity of suitable weak solutions of the Navier-Stokes
equations}, Comm. Pure Appl. Math. {\bf 35} (1982), 771-831.



\bibitem{Farwig_Sohr} {\sc R.~Farwig, H.~Sohr}, {\it The
stationary and nonstationary Stokes system in exterior domains
with nonzero divergence and nonzero boundary data,} Math. Meth.
Appl. Sci. {\bf 17} (1994), 269-291.

\bibitem{LSU}
{\sc O.A.~Ladyzhenskaya, V.A.~Solonnikov, N.N.~Uraltseva}, Linear
and quasilinear equations of parabolic type.
 Translations of Mathematical Monographs, Vol. 23 American Mathematical
 Society, Providence, R.I., 1967.

\bibitem{Seregin_ZNS_271}
{\sc G.A.~Seregin}, {\it  Some estimates near the boundary for
solutions to the non-stationary linearized Navier-Stokes
equations}, Zapiski Nauchnyh Seminarov POMI {\bf 271} (2000),
204-223.

\bibitem{Seregin_JMFM}
{\sc G.A.~Seregin}, Local regularity of suitable weak solutions to
the Navier-Stokes equations near the boundary, Journal of
Mathematical Fluid Mechanics  {\bf 4} (2002) no.1, 1-29.

\bibitem{Seregin_Handbook}
{\sc G.A.~Seregin}, {\it ``Local regularity theory of the Navier-Stokes equations''}, Handbook of Mathematical Fluid Dynamics, Volume 4 (2007), 159-200.

\bibitem{Solonnikov_PMA_21}
{\sc V.A. Solonnikov}, {\it Estimates in $L_p$ of solutions to the initial-boundary value problem for the generalized Stokes system in a bounded domain},
Problems of Math. Analysis {\bf 21} (2000), 211-263.

\bibitem{Solonnikov_ZNS_288}
{\sc V.A. Solonnikov}, {\it Estimates of solutions of the Stokes
equations in Sobolev spaces with a mixed norm}, Zapiski Nauchnyh
Seminarov POMI {\bf 288} (2002), 204-231.

\bibitem{Solonnikov_Uspekhi}
{\sc V.A. Solonnikov}, {\it On the estimates of solutions of
nonstationary Stokes problem in anisotropic Sobolev spaces
and on the estimate of resolvent of the Stokes problem},
Uspekhi Matematicheskih Nauk, {\bf 58} (2003) no.2 (350), 123-156.

\bibitem{Sol_Ur}
{\sc D.K.~Faddeev, B.Z.~Vulich, V.A.~Solonnikov, N.N.~Uraltseva,} {\it Izbrannye glavy analiza i vyschei algebry}, Izdatel'stvo LGU, 1981 (in Russian).



\end{thebibliography}
\end{document}